\documentclass[10pt,a4paper]{article}

\usepackage{color}
\usepackage{ifpdf}
\ifpdf
 \usepackage[pdftex]{graphicx}
\else
 \usepackage[dvips]{graphicx}
\fi

\usepackage{hyperref}

\usepackage[T1]{fontenc}
\usepackage{a4,amssymb,amsthm,amsmath}
\usepackage{lmodern}
\usepackage[all]{xy}
\def\à{\`a}
\def\è{\`e}
\def\ä{\"a}

\newcommand{\g}{\mathfrak{h}}
\newcommand{\R}{\mathbb{R}}
\newcommand{\C}{\mathbb{C}}

\newcommand{\X}{\mathcal{X}}
\newcommand{\Y}{\mathcal{Y}}

\newcommand{\T}{\mathbb{T}}

\newcommand{\V}{\mathcal{V}}
\newcommand{\NN}{\mathcal{N}}
\newcommand{\QQ}{\mathcal{Q}}
\newcommand{\HH}{\mathcal{H}}

\newcommand{\Sp}{\mathbb{S}}
\newcommand{\Z}{\mathbb{Z}}
\newcommand{\ZZ}{\mathcal{Z}}

\newcommand{\OO}{\mathcal{O}}
\newcommand{\RR}{\mathcal{R}}
\newcommand{\TT}{\mathcal{T}}
\newcommand{\II}{\mathcal{I}}
\newcommand{\KK}{\mathcal{K}}

\newcommand{\J}{\mathcal{J}}
\newcommand{\K}{\mathcal{K}}
\newcommand{\JJ}{\mathbb{J}}
\newcommand{\Q}{\mathbb{H}}

\newcommand{\lra}{\longrightarrow}
\newcommand{\lms}{\longmapsto}
\newcommand{\bw}{\bigwedge}
\newcommand{\w}{\wedge}
\newcommand{\wh}{\widehat}

\theoremstyle{definition}

\author{Guillaume~Deschamps}
 \title{   Espace des twisteurs d'une variété quaternionique Kähler généralisée.}
 \date{\today}

\begin{document}

\maketitle

\begin{center}
{\bf Résumé}
\end{center}
Munir une variété $M$ de dimension $4n$, d'une structure presque
quaternionique $Q$ revient précisément à lui associer un fibré des
twisteurs $Z(Q)\lra M$. Lorsque $Q$ est stable par une connexion sans
torsion, on peut munir $Z(Q)$ d'une structure presque complexe
$\JJ$. Dans le cas $n=1$, les travaux d'Atiyah, Hitchin et Singer
\cite{AHS78} ont permis de relier l'intégrabilité de $\JJ$ à la
géométrie de la variété $(M,Q)$. Pour $n>1$, Salamon \cite{Sal,
Sal1} a montré que la structure presque complexe $\JJ$ sur
$Z(Q)$ est toujours intégrable. Pantilie \cite{Pan} a remarqué qu'on pouvait étendre ces résultats à la géométrie complexe généralisée. Ainsi il définit ce qu'est une variété presque
quaternionique généralisée $(M,\QQ)$ et lui associe un $\Sp^2$-fibré  $\ZZ(\QQ)\lra M$. Comme dans le cas
classique, lorsque $\QQ$ est stable par une connexion sans torsion
(au sens des connexions généralisées), il munit $\ZZ(\QQ)$
d'une structure presque complexe généralisée $\mathbb J$ mais ne donne pas de critère d'intégrabilité.  Le but de cette article est précisément de donner  
 un critère d'intégrabilité pour $\mathbb J$. Nous étudierons ensuite plus particulièrement le cas
où $(M,g,\QQ)$ est une  variété quaternionique Kähler
généralisée et verrons qu'alors $\JJ$ est automatiquement  intégrable, sous réserve que 
$n>1$. Nous illustrerons ces résultats en donnant plusieurs
exemples.

\maketitle \setcounter{tocdepth}{4}

\begin{center}
{\bf Abstract}
\end{center}
To give an almost quaternionic structure on a 4n-manifold $M$ is
equivalent to give its bundle of twistors $Z(Q)\lra M$. When $Q$
 is invariant under a torsion free
connection, $Z(Q) $ can be provided with an almost complex
structure $ \JJ $. In the case $ n = 1 $  Atiyah, Hitchin and
Singer \cite{AHS78} have related the integrability of $ \JJ $ to
the geometry of  $ (M, Q) $. For $ n> 1 $ Salamon \cite{Sal, Sal1}
showed that the almost complex structure $ \JJ $ on $ Z (Q) $
is always integrable. Pantilie \cite{Pan} notice that 
these results can be extended to the generalized complex geometry: he defines the concept of almost generalized quaternionic manifold
$ (M, \QQ ) $ and its twistor space $\ZZ(\QQ)$, which is a $\Sp^2$-fiber bundle over $M$. As in the usual case, when $\QQ$ is invariant under a
generalized torsion free connection, he shows that  $ \ZZ(\QQ) $ comes
with an almost generalized complex structure $\JJ$ but doesn't give any integrability criterion. The purpose of this article is precisely to give an integrability criterion for $\JJ$. In particular, when $ (M, g,\QQ) $ is a
generalized quaternionic Kähler manifold, we will show that $\JJ$ is
always integrable as soon as $n>1$. We illustrate this work by
giving several examples.

\newpage
\setcounter{tocdepth}{4}
\tableofcontents

\section{Introduction}
 Introduite par Penrose \cite{Pen76}, la théorie des twisteurs a
pris tout son essor avec l'article d'Atiyah, Hitchin et Singer
\cite{AHS78} qui, à toute 4-variété riemannienne orientée $(M,g)$
associe un fibré $Z(M,g)\lra M$ en sphère $\Sp^2$. L'espace
$Z(M,g)$ est alors muni d'une structure presque complexe
naturelle $\JJ$ dont l'intégrabilité dépend de la métrique $g$ sur
$M$. En dimension plus grande et pour une variété presque
quaternionique $(M,Q)$, on peut toujours définir un fibré des
twisteurs $Z(Q)\lra M$ qui est encore un fibré en sphère
$\Sp^2$. Lorsque $Q$ est  stable par une connexion sans torsion on
peut munir $Z(Q)$ d'une structure presque complexe $\JJ$. Dans
ce cas Salamon \cite{Sal,Sal1} a montré que l'intégrabilité de
$\JJ$ est automatique.

C'est alors qu'a été introduit par Hitchin \cite{Hit} le concept
de structures complexes généralisées, dans le but d'unifier les
notions de structures complexes et de structures symplectiques. En
utilisant cette théorie, Pantilie \cite{Pan} a pu définir de façon
très naturelle la notion de variétés presque quaternioniques
généralisées $(M,\mathcal Q)$ et d'espace de twisteurs associé
$\ZZ(\mathcal Q)$. Contrairement à certaines généralisations \cite{Bre, Des2, DM1,DM2, GS},
 ici l'espace des twisteurs reste un fibré en sphères $\Sp^2$.
 Lorsque $\QQ$ est
stable par une connexion sans "torsion généralisée" (cf. section
2.3 pour des définitions précises), la variété $\ZZ(\mathcal Q)$ admet
une structure presque complexe généralisée, que nous noterons
encore $\JJ$.

 Dans son article, Pantilie donne deux exemples
d'espace de twisteurs pour lesquelles il vérifie que la
structure presque complexe généralisée $\JJ$ est intégrable mais ne donne pas de critère d'intégrabilité. Pour palier à cela, on se
propose dans cet article d'établir un théorème d'intégrabilité similaire à ceux
d'Atiyah, Hitchin, Singer et de Salamon c'est-à-dire qui exprime
l'intégrabilité de la structure presque complexe généralisée $\JJ$
sur $\ZZ(\QQ)$ en fonction de la géométrie de la variété
$(M,\QQ)$. Pour cela nous
 rappellerons dans la partie 2, les différentes
 constructions d'espaces de twisteurs.
 On dira en particulier ce qu'est une structure presque quaternionique généralisée
$\QQ$ sur une variété  $M$; ce qu'est l'espace des twisteurs
$\ZZ(\QQ)\lra M$ associé à $(M,\QQ)$ et enfin lorsque $\QQ$ est
une structure quaternionique généralisée, nous rappellerons la construction de la structure presque complexe généralisée $\JJ$ sur $\ZZ(\QQ)$.
 Dans la partie 3 nous énoncerons le théorème d'intégrabilité de
 $\JJ$, nous verrons qu'il dépend de la courbure de la connexion
 qui définit $\JJ$. Dans la partie 4, nous appliquerons ce
 résultat dans le cas particulier des variétés quaternioniques
 Kähler généralisées :

  \quad\\
  {\bf Théorème.} Soit $(M,g,\QQ)$ une 4n-variété quaternionique
  Kähler généralisée avec $n>1$. La structure presque complexe
  généralisée $\JJ$ sur son espace de twisteurs $\ZZ(\QQ)$ est
  toujours intégrable.

\quad\\
Pour les 4-variétés, l'intégrabilité de $\JJ$ dépend de la
courbure de la métrique $g$ (cf théorème 3 et 4). 
Enfin dans la partie 5,  nous illustrerons ces résultats  par des exemples.

\section{Rappels et définitions}
Soit $E\lra M$ un fibré au-dessus de $M$,  nous noterons $\Gamma(E)$ l'ensemble des sections lisses.
 \subsection{Structure quaternionique.} Une
structure {\it presque complexe} sur une variété $M$ est un
endomorphisme du fibré tangent $J\in \Gamma(End(TM))$ qui satisfait
$J^2=-Id$. Cela revient à munir les fibres du tangent $TM\lra M$
d'une structure de $\C$-espace vectoriel et force donc la
dimension de $M$ à être paire. Si la variété $M$ est complexe elle
est munie d'une structure presque complexe canonique.
Réciproquement, on dira qu'une structure presque complexe $J$ sur
$M$ est {\it intégrable} si elle provient d'un atlas holomorphe.
Le théorème de Newlander-Nirenberg nous donne un critère
d'intégrabilité.

\quad\\
{\bf Théorème \cite{NN}.} Soit $(M,J)$ une variété presque
complexe. Les conditions suivantes sont équivalentes :

\begin{enumerate}
\item La structure presque complexe $J$  est intégrable.

\item  Le tenseur de Nijenhuis défini par :
$$
 N(X,Y)=[JX,JY]-J[JX,Y]-J[X,JY]-[X,Y]
$$
 est nul pour toutes sections $X,Y\in\Gamma(TM)$.

\item Le sous-espace espace propre $T^{1,0}$ de $J$ dans
$TM\otimes\C$ associé à la valeur propre $i$,  est stable par
crochet de Lie (i.e. $T^{1,0}$ involutif).

\end{enumerate}
 \quad\\
 Soit $(M,g)$ une variété riemannienne munie d'une structure
presque complexe $J$. Nous dirons que $J$ est compatible avec $g$
si $J$ est un endomorphisme orthogonal pour $g$ et nous noterons
$O_g(TM)$ l'ensemble des endomorphismes orthogonaux pour $g$. Nous
noterons plus généralement $End(TM)$ l'ensemble des endomorphismes
de $TM$. Lorsque $J$ est compatible avec $g$, on peut définir la
2-forme fondamentale 
$$
w(X,Y)=g(JX,Y)\quad\forall X,Y\in \Gamma(TM).
$$
Nous dirons que $(M,g,J)$ est une variété {\it kählérienne} si $J$
est intégrable et si $w$ est fermée.

Supposons maintenant que la variété $M$ admette deux structures
presque complexes $I$ et $J$, qui anti-commutent. Si on pose
$K=IJ$ alors $K$ est une nouvelle structure presque complexe sur
$M$ et plus généralement pour tout point $(a,b,c)\in\Sp^2$
l'endomorphisme $aI+bJ+cK$ définit aussi une structure presque
complexe. C'est pourquoi nous dirons que la variété $(M,I,J,K)$
est {\it presque hypercomplexe}. Elle sera dite hypercomplexe dès
que $I$ et $J$ sont intégrables. Dans ce cas toute la
$\Sp^2$-famille de structures presque complexes est intégrable.
Par ailleurs, nous dirons que la variété $(M,g,I,J)$ est {\it
hyperkählérienne} lorsque les structures complexes $I$ et $J$ sont
kählériennes.  Là aussi, la terminologie se justifie car pour tout
point $(a,b,c)\in \Sp^2$, la structure presque complexe $aI+bJ+cK$
est kählérienne. Dans tous ces cas, la dimension de $M$ est
nécessairement un multiple de quatre.

 \quad\\
 {\bf Définition.}   Une structure {\it presque
quaternionique} sur une variété $M$ est un sous-fibré de rang
trois $Q\subset End(TM)$ localement engendré par une structure
presque hypercomplexe.

\quad\\
Une connexion sur $TM$ induit une connexion sur $End(TM)$ en
posant, pour tout $\psi\in \Gamma(End(TM))$ et pour tous $X,Y\in \Gamma(TM)$ :
$$
(\nabla_X\psi)(Y)=\nabla_X \psi(Y)-\psi(\nabla_X Y).
$$
Lorsqu'il existe une connexion sans torsion sur $TM$ qui stabilise
$Q$ au sens où $\nabla Q\subset Q$,  on parle d'une structure {\it
quaternionique}. Une structure presque quaternionique sur $M$
équivaut à se donner une $GL(n,\mathbb{H})Sp(1)$-structure sur
$M$. Et une connexion sans torsion sur $TM$ qui stabilise $Q$
revient à se donner une $GL(n,\mathbb{H})Sp(1)$-connexion
  sans torsion \cite{Bes87}.
  Lorsque   $(M,g,Q)$ est une variété  presque quaternionique munie
d'une  métrique riemannienne $g$ telle que  $Q$ soit localement
engendré par une structure presque hypercomplexe, où $I,J$ et $K$
sont compatibles avec $g$; on dit que $(M,g,Q)$ est une structure
{\it presque quaternionique hermitienne}.  Cela correspond à se
donner une $Sp(n)Sp(1)$-structure sur $M$. Si de plus $Q$  est
stable par la connexion de Levi Civita on dit que $(M,g,Q)$ est
une variété {\it quaternionique Kähler}. Une variété
quaternionique Kähler est donc précisément une variété
riemannienne dont le groupe d'holonomie est inclus dans
$Sp(n)Sp(1)$. Lorsque $n>1$, Alekseevskii et Berger ont montré que
cela imposait à la métrique $g$ d'être Einstein \cite{Ale, Ber}.
En particulier si on regarde le tenseur de courbure $R$ associé à
la métrique $g$ : $R(X,Y)=\nabla_X\nabla_Y-\nabla_Y\nabla_X-\nabla_{[X,Y]}$; vu
comme endomorphisme auto-adjoint du  fibré extérieur $\bw^{2}T M$
on a 
$$
 R\vert_Q=\lambda Id\vert_Q,
$$
où $\lambda$ est un multiple positif de la courbure scalaire de
$(M,g)$ et où on identifie les éléments  $u$ de $Q$ aux éléments
$\phi(u)$ de $\bw^{2} TM $ via
$$g\Big(\phi(u),X\w Y\Big)=g(uX,Y)
\quad\forall X,Y\in TM.$$ En dimension 4 une structure presque
quaternionique revient à se donner une structure conforme orientée
et les structures presque quaternioniques sont automatiquement
stables par toute connexion.
 Plus précisément, si $(M,g)$ est une 4-variété riemannienne,
  l'opérateur étoile de Hodge induit la décomposition $\bw^{2}T M=\bw^+\oplus\bw^-$ et les
seules structures presque quaternioniques hermitiennes de $TM$
sont précisément $\bw^+$ et $\bw^-$. Notons enfin que dans cette
base, le tenseur de courbure $R:\bw^2T M\lra\bw^2T M$ se décompose
en blocs \cite{ST, Bes87} :
$$
R=\left[\begin{array}{cc}\mathcal W^++\frac{s}{12}Id&\mathcal B\\
\mathcal B^\star&\mathcal W^-+\frac{s}{12}Id\end{array}\right].
$$
L'opérateur $\mathcal W=\mathcal W^++\mathcal W^-$ est l'opérateur
de Weyl, $s$ la courbure scalaire, $\mathcal B$ le tenseur de
Ricci sans trace, $\mathcal B^\star$ son adjoint  et $Id$ la
matrice identité.

\quad\\
Lorsque $(M,Q)$ est une variété presque quaternionique de
dimension $4n$ avec $n\geq1$, on peut définir le fibré des
twisteurs $\pi:Z(Q)\lra M$ comme le fibré des structures presque
complexes sur $M$ appartenant à $Q$. C'est un fibré de fibres
$\Sp^2$ et de groupe structural $SO(3)$. Les fibres admettent donc
naturellement une structure complexe. Une connexion
 qui stabilise $Q$ fournit alors une décomposition de l'espace tangent
de $Z(Q)$ :
$$
TZ(Q)=\mathcal H\oplus \mathcal V
$$
 en une distribution horizontale $\mathcal H$ et une distribution verticale
  $\mathcal V=ker\,d\pi$ (i.e. tangente aux fibres).
  En un point $p$ de $Z(Q)$, comme
$\mathcal H_{p}$ est isomorphe à $T_{\pi(p)}M$ via $d\pi$, la
distribution horizontale hérite naturellement de la structure
presque complexe induite par $p$. La somme de cette structure
presque complexe et de celle sur les fibres munit $Z(Q)$ d'une
structure presque complexe naturelle notée $\JJ$. En dimension
quatre, un critère d'intégrabilité de $\JJ$  a été donné par
Atiyah, Hitchin et Singer.

\quad\\
{\bf Théorème \cite{AHS78}.} Soit  $(M,g)$ une 4-variété
riemannienne orientée. La structure presque complexe $\JJ$ sur
$Z(\bw^+)$ associée à la connexion de Levi-Civita, est
intégrable si et seulement si $g$ est anti-autoduale (i.e.
$\mathcal W^+=0$).

\quad\\
 En dimension plus grande, l'intégrabilité de $\JJ$ a été démontré par
 Salamon.

\quad\\
{\bf Théorème \cite{Sal, Sal1}.} Soit $n> 1$ et $(M,Q)$ une
4n-variété quaternionique alors la structure presque complexe
$\JJ$ sur $Z(Q)$ est toujours intégrable.

\quad\\
{\bf Remarque.} Pour une 4-variété riemannienne orientée $(M,g)$,
la structure presque complexe $\JJ$ sur $Z(\bw^+)$ associé à la
connexion de Levi-Civita ne dépend que de la classe conforme de
$g$. Pour que ces deux théorèmes coïncident, certains auteurs
définissent donc les 4-variétés quaternioniques comme des
4-variétés munies de structures conformes orientées
anti-autoduales. Pour d'autres les variétés quaternioniques sont
de dimension $4n>4$.

\subsection{Structure complexe généralisée.} Soit $M$ une
variété de dimension $2n$. En géométrie généralisée on étudie non
pas le fibré  tangent
  de M  mais  la somme  du fibré tangent et du fibré cotangent que
   nous noterons $\T M=TM\oplus T^\star M$. Sur $\T M$ il y a
une pseudo-métrique naturelle de signature $(2n,2n)$ définie par :
$$
<X+\xi,Y+\eta>=\frac{1}{2}\Big(\xi(Y)+\eta(X)\Big) \quad\forall
X,Y\in TM\textrm{ et }\forall\xi,\eta\in  T^\star M.
$$
Sur $\T M$ on a également le crochet de Courant \cite{Cou}, défini
pour toutes sections $X+\xi,Y+\eta\in \Gamma(\T M)$ par
$$
[X+\xi,Y+\eta]=[X,Y]+\mathcal L_X\eta-\mathcal
L_Y\xi-\frac{1}{2}d(i_X\eta-i_Y\xi)
$$
où $[X,Y]$ est le crochet de Lie, $\mathcal L_X$ la dérivée de Lie
et $i_X$ le produit intérieur (i.e. la contraction par le premier
argument).

\quad\\
{\bf Définition.}
 Une
structure {\it presque complexe généralisée} \cite{Hit, Gua2} sur
$M$ est la donnée d'un endomorphisme du fibré tangent généralisée $\J\in\Gamma(End(\T M))$ telle que
\begin{enumerate}
\item[i)] $\J^2=-Id$

\item[ii)] $\J$ préserve la pseudo-métrique $<.,.>$.

\end{enumerate}

\quad\\
 {\bf Théorème-Définition \cite{Gua2}.} Une structure presque complexe généralisée
$\J$ sur $M$ est intégrable si et seulement si l'une des deux
conditions (équivalentes) suivantes est vérifiée:

\begin{enumerate}
\item  Le tenseur de Nijenhuis défini par 
$$
\NN(\mathcal X, \mathcal Y)=[\J \mathcal X,\J \mathcal Y]-\J[\J
\mathcal X,\mathcal Y]-\J[ \mathcal X,\J \mathcal Y]-[\mathcal
X,\mathcal Y]
$$
est nul pour toutes sections $\mathcal X, \mathcal Y$ de $\T M$.

\item Le sous-espace propre $\T ^{1,0}$ de $\J$ dans $\T
M\otimes\C$ associé à la valeur propre $i$,  est stable par
crochet de Courant.
\end{enumerate}
 \quad\\
 On note $pr_1 : (TM\oplus T^\star M)\otimes\C\lra TM\otimes\C$
la première projection. La codimension dans $TM\otimes\C$ de
$pr_1(\T^{1,0})$ est un invariant de la structure presque complexe
appelé le {\it type} de $\J$.

 Comme le montre les exemples suivants, la notion de structure complexe généralisée
  regroupe sous
un même formalisme les notions de structures complexes et de
structures symplectiques.

\quad\\
{\bf Exemple a.}  Une structure presque complexe $J$ sur $M$
définit la structure presque complexe généralisée $\mathcal
J_J=\left(\begin{array}{cc} J&0\\0&-J^\star\end{array}\right)$ où
$J^\star$ est l'adjoint de $J$. De plus $\mathcal J_J$ est
intégrable si et seulement si $J$ l'est. Le type de $\mathcal J_J$
est constant égale à $n$.

\quad\\
{\bf Exemple b.} De même une structure presque symplectique $w$
sur $M$ (i.e. une 2-forme non dégénérée),  définit la structure
presque complexe généralisée $\mathcal J_w=\left(\begin{array}{cc}
0&-w^{-1}\\w&0\end{array}\right)$. La structure  $\mathcal J_w$
est intégrable si et seulement si $w$ l'est (i.e. $w$ fermée). Le
type de $\mathcal J_w$ est $0$.

\quad\\
{\bf Exemple c.} Toute 2-forme $B$ sur $M$ définit  l'application
orthogonale suivante:
$$
\begin{array}{lccc}
e^B :&TM\oplus T^\star M&\lra&TM\oplus T^\star M\\
&X+\xi&\lms&X+\xi+i_X B.
\end{array}
$$
 Si $\J$ est une
structure presque complexe généralisée sur $M$ alors $e^{-B}\J
e^B$ aussi et de même type que $\J$. De plus lorsque $B$ est
fermée, $\J$ est intégrable si et seulement sa $B$-transformation
$e^{-B}\J e^B$ l'est.

\quad\\ Une structure complexe sur une $2n$-variété est par
définition, localement difféomorphe à $\C^n$. D'autre part, le
théorème de Darboux nous dit qu'une structure symplectique sur une
$2n$-variété est localement difféomorphe à $(\R^{2n},w_0)$ où :
$$
w_0=dx_1\w dx_2+\ldots,+dx_{2n-1}\w dx_{2n}.
$$
De même en un point régulier (i.e. où le type est constant égale à
$k$ dans un voisinage du point) une structure complexe généralisée
est, à B-transformation près,  localement difféomorphe au produit
de $\C^k$ et de $(\R^{2n-2k},w_0)$ \cite{Gua2}.

\subsection{Structure quaternionique généralisée.}
Comme l'a remarqué Pantilie \cite{Pan}, les définitions données dans la partie 2.1 s'étendent naturellement
dans le cadre de la géométrie complexe généralisée. Ainsi, une variété $M$ est {\it presque hypercomplexe
généralisée} s'il existe une paire de structures presque complexes
généralisées $\II, \J$ qui anti-commutent. Là encore, si on pose
 $\KK=\II\J$ alors $\KK$ définit une structure presque complexe généralisée,
  tout comme l'endomorphisme $a\II+b\J+c\KK$ associé
 au point $(a,b,c)\in\Sp^2$.

 \quad\\
 {\bf Définition.}   Une structure {\it presque
quaternionique généralisée} sur une variété $M$ est un sous-fibré
de rang trois $\QQ\subset End(\T M)$ localement engendré par une
structure presque hypercomplexe généralisée.

\quad\\ Nous noterons $\overrightarrow \X\in T M$ la partie
vectorielle de $\X\in \T M$, c'est-à-dire la projection sur $TM$
parallèlement à $T^\star M$:
$$\begin{array}{lccc}
\overrightarrow{}:&\T M&\lra& T M\\
&X+\xi&\lms&X.
\end{array}
$$
 Une connexion (classique) $\nabla$ sur $\T M$ s'étend
naturellement  en une connexion généralisée   \cite{Gua3} sur $\T M$ en posant :
$$
\nabla_\X\Y:=\nabla_{\overrightarrow \X}\Y
\quad \forall\X,\Y\in \Gamma(\T M).
$$
 La torsion d'une
connexion généralisée sur $\T M$ est définie pour toutes sections
$\X_1,\X_2,\X_3$ de  $\T M$ par \cite{Gua3}:
$$
\begin{array}{ccc}
\TT(\X_1,\X_2,\X_3)&=&<\nabla_{\X_1}\X_2-\nabla_{\X_2}\X_1-[\X_1,\X_2],\X_3>\\
&& +\displaystyle
\frac{1}{2}\Big(<\nabla_{\X_3}\X_1,\X_2>-<\nabla_{\X_3}\X_2,\X_1>\Big).
\end{array}$$
On appellera  torsion généralisée de $\nabla$ l'opérateur $\TT$.

\quad\\
{\bf Définition.} S'il existe une connexion (classique) $\nabla$ sur $\T M$
compatible avec la pseudo-métrique, qui stabilise $\QQ$ (i.e.
$\nabla \QQ\subset \QQ$) et qui est sans torsion généralisée; on dit
que $(M,\QQ,\nabla)$ est une variété {\it quaternionique
généralisée}.

\quad\\
Une métrique riemannienne $g$ sur $M$ se prolonge en une métrique
sur $\T M$. Une structure presque complexe généralisée est dite
hermitienne ou compatible avec $g$ si elle agit comme un endomorphisme orthogonal
pour $g$.

\quad\\
 {\bf Définition.} Une variété riemannienne $(M,g,\QQ)$ est dite
 {\it quaternionique Kähler généralisée} si $\QQ$ est une structure presque
quaternionique hermitienne généralisée, stable par la connexion de
Levi-Civita.

\quad\\
Comme dans le cas classique, on peut donner une interprétation de
ces définitions en termes de $G$-structure. Pour cela, notons
$\R^{4n\star}$ le dual de $\R^{4n}$ et $\mathcal O(\T M)$ le
sous-groupe de $End(\T M)$ des isométries pour la pseudo-métrique.
Le groupe $GL(2n,\Q)$ agit naturellement à gauche sur
$\R^{4n}\oplus\R^{4n\star}$. On notera $\mathcal S p(2n)$ le
sous-groupe de $GL(2n,\Q)$ des isométries pour la
pseudo-métriques:
$$
\mathcal S p(2n)=GL(2n,\Q)\cap\mathcal
O\Big(\R^{4n}\oplus\R^{4n\star}\Big).
$$
 Le groupe $Sp(1)$ des quaternions de
norme 1 agit naturellement à droite sur
$\R^{4n}\oplus\R^{4n\star}$. On notera $G$ le groupe produit
$\mathcal S p(2n)Sp(1)$.
  Une structure quaternionique généralisée
$\QQ$ sur $M$  équivaut à se donner un G-fibré principal $\mathcal
P\lra M$ qui est un sous-fibré de $\mathcal O(\T M)\lra M$, ainsi
qu'une $G$-connexion sur $\mathcal P$  sans torsion généralisée.

\subsection{Espace des twisteurs.}
 Lorsque $(M,\QQ)$ est une  variété presque quaternionique généralisée,
 on peut lui associer tout naturellement un fibré des
twisteurs $\pi : \ZZ(\QQ)\lra M$. C'est le fibré des structures
presque complexes généralisées sur $M$ appartenant à $\QQ$. C'est
un fibré en sphère $\Sp^2$ et de groupe structural $SO(3)$. C'est
aussi le fibré associé au fibré principal $\mathcal P\lra M$ par
l'action par conjugaison à droite de $G$ sur $\Sp^2$ (c.f. section
3.3).

\quad\\
 Comme $\ZZ(\QQ)$ est un fibré associé, une $G$-connexion
  sur $\mathcal P\lra M$ permet de décomposer
l'espace tangent $T\ZZ(\QQ)$ en la somme directe d'une
distribution horizontale $\HH$ et de la distribution verticale
$\V=ker\,d\pi$. On identifie le dual de $\HH$ (resp. de $\V$) noté
$\HH^\star$ (resp. $\V^\star$) aux formes linéaires sur
$T\ZZ(\QQ)$ nulles sur $\V$ (resp. sur $\HH$). Comme dans le cas
classique, on peut munir $\ZZ(\QQ)$ d'une structure presque
complexe généralisée naturelle. En effet en un point $p$ de
$\ZZ(\QQ)$ comme $\HH_p\oplus\HH^\star_p$ est isomorphe à $\T
_{\pi(p)} M$ via $d\pi\oplus d\pi^\star$, il hérite de la
structure presque complexe généralisée induite par $p$. La somme
de cette structure presque complexe généralisée et de la structure
complexe sur les fibres munit $\ZZ(\QQ)$ d'une structure presque
complexe généralisée encore notée $\JJ$.

\quad\\
 La question  est de savoir sous quelle condition
 la structure presque complexe généralisée $\JJ$ sur $\ZZ(\QQ)$ est intégrable? Avant de répondre à cette
question,  rappelons les  exemples proposés par Pantilie \cite{Pan}.

\quad\\
{\bf Exemple 1.} Comme pour l'exemple a, une variété
quaternionique (Kähler) est naturellement quaternionique (Kähler)
généralisée.

\quad\\
{\bf Exemple 2 \cite{Pan}.} Soit $(M,g,I,J,K)$ une variété
hyperkählérienne. Cette variété étant hyperkählérienne elle est
aussi, quaternionique Kähler généralisée. Mais on peut aussi
construire une structure "tordue". Ainsi, notons
$\J_I=\left(\begin{array}{ll}I&0\\0&I\end{array}\right)$ et
$\J_w=\left(\begin{array}{ll}0&J\\J&0\end{array}\right)$ les
structures complexes généralisées associées à la structure
complexe $I$ et à la structure symplectique $w=g(J.,.)$. Comme $M$
est hyperkählérienne on a
$\J_I\J_w=-\J_w\J_I=\left(\begin{array}{ll}0&K\\K&0\end{array}\right)$.
On note $\J_{Iw}$ cette dernière structure complexe généralisée.
La distribution  $\QQ=\{a\J_I+b\J_w+c\J_{Iw}/(a,b,c)\in\R^3\}$
définit une structure quaternionique Kähler généralisée sur
$(M,g)$. Pour une telle variété, il est facile de vérifier que la
structure presque complexe généralisée $\JJ$ sur $\ZZ(\QQ)$ est
intégrable.

\quad\\
{\bf Exemple 3 \cite{Pan}.} Soit $(N,c,\nabla)$ une 3-variété
Einstein-Weyl et $\varphi :(M,g)\lra (N,c,\nabla)$ son "heaven
space".  La 4-variété $M$ est alors parallélisable et la métrique
$g$ est anti-autoduale, Einstein à courbure scalaire non nul.
Comme toute 4-variété riemannienne orientée, $(M,g)$ est
naturellement quaternionique Kähler et donc quaternionique Kähler
généralisée. Mais là encore, on peut munir $(M,g)$ d'une structure
"tordue". En effet l'application $\varphi : M\lra N$ induit une
décomposition $TM\simeq\mathcal H\oplus \mathcal V^\star$ avec
$\mathcal V=ker d\varphi$ et $\mathcal H=\mathcal V^\perp$, de
sorte que sous cet isomorphisme, la structure quaternionique
Kähler naturelle devient une structure quaternionique généralisée
$\QQ$. La structure presque complexe généralisée $\JJ$ sur
$\ZZ(\QQ)$ est alors intégrable.

\section{Espace des twisteurs et critère d'intégrabilité}
\subsection{Opérateur de courbure.} Soit $(M,\QQ,\nabla)$ une variété
quaternionique
 généralisée et $\pi:\ZZ(\QQ)\lra M$ l'espace des twisteurs associé.
  Soit $\mathcal U$ un petit ouvert de $M$ sur lequel on a une
  trivialisation  $\pi^{-1}({\mathcal U})\simeq \mathcal U\times
  \Sp^2$ et $(m,u)\in \mathcal U\times
  \Sp^2$ un système de coordonnées locales.
Une connexion généralisée admet un opérateur de courbure $\RR$ défini sur $\T M$ par \cite{Gua3}:
$$
\RR(\X,\Y)\mathcal Z=[\nabla_\X,\nabla_\Y]\mathcal Z
-\nabla_{[\X,\Y]}\mathcal Z, \quad \forall \X,\Y, \mathcal Z\in
\Gamma(\T M).
$$
Par anti-symétrie en $\X$ et $\Y$ on écrira parfois $\RR(\X\w\Y)$
plutôt que $\RR(\X,\Y)$ et pour alléger l'écriture, si $\X,\Y\in
\Gamma\Big(\T \ZZ(\QQ)\Big)$ on écrira $\RR(\X,\Y)$
plutôt que $\RR(\pi_\star\X,\pi_\star\Y)$.

\quad\\
{\bf Remarque.} Les connexions que nous considérons dans cette article vérifient  $\nabla_\xi=0$ pour tout $\xi\in\Gamma(T^\star M)$ si bien que  $\RR$ est linéaire en $\X$ et $\Y$.

\quad\\
La connexion $\nabla$ agit sur $End(\T M)$; pour tous
$\X,\Y\in\Gamma(\T M)$ on peut donc faire agir
l'opérateur de courbure $\RR(\X,\Y)$ sur $End(\T M)$. Comme $\nabla$
stabilise $\QQ$; pour tout $u\in \Gamma(\QQ)$ on a toujours $\RR(\X,\Y)u\in\Gamma(\QQ)$. Mais l'opérateur $\RR(\X,\Y)$ est par définition un
élément de $End(\T M)$. On peut donc aussi regarder le commutateur d'un élément $u\in\QQ$ avec $\RR(\X,\Y)$ :
$$[u,\RR(\X,\Y)]=u \RR(\X,\Y)-\RR(\X,\Y)u.
$$
 On vérifie qu'on a l'égalité:
$$[u,\RR(\X,\Y)]=-\RR(\X,\Y)u,\quad\forall u\in\Gamma(\QQ).
$$
Le commutateur $[u,\RR(\X,\Y)]$ est donc  un élément de $\QQ$, et plus précisément, on peut voir que c'est un champ de vecteurs vertical sur $\ZZ(\QQ)$.
\subsection{\'Enoncé.}
 Le résultat central de cet article est le suivant:

\quad\\
{\bf Théorème 1.}  Soit $(M,\QQ,\nabla)$ une 4n-variété  munie
d'une structure quaternionique
 généralisée avec $n\geq1$. La structure presque complexe généralisée $\JJ$
sur $\ZZ(\QQ)$  est intégrable si et seulement si pour tous
champs $\X,\Y\in \Gamma(\HH\oplus\HH^\star)$ et pour tout point
$(m,u)\in\pi^{-1}(\mathcal U)$ on a :
$$\left[u,\RR\Big(\X\wedge\Y-u\X\wedge  u\Y\Big)
+u\RR\Big(u\X\wedge\Y+\X\wedge u\Y\Big)\right]=0.$$

\quad\\
 {\bf Notation.} Soit $\X\in \T \ZZ(\QQ)$,  on note $\X^{1,0}\in \T^{1,0}$ (resp.
 $\X^{0,1}\in\T^{0,1}$) la partie $(1,0)$ (resp. $(0,1)$) du champ
 $\X$,
 c'est-à-dire la projection sur le sous-espace propre de $\JJ$
 associé à la valeur propre $i$ (resp. $-i$) :
$$
\X^{1,0}=\frac{\X-i\JJ\X}{2}\textrm{ et }
\X^{0,1}=\frac{\X+i\JJ\X}{2}.
$$
Avec cette notation, on peut reformuler le théorème précédent de
manière plus intrinsèque.

\quad\\
{\bf Théorème 1 bis.}  Soit $(M,\QQ,\nabla)$ une 4n-variété  munie
d'une structure quaternionique
 généralisée avec $n\geq1$. La structure presque complexe généralisée $\JJ$
sur $\ZZ(\QQ)$  est intégrable si et seulement si :
$$\Big(\RR(\X^{1,0},\Y^{1,0})\mathcal Z^{1,0}\Big)^{0,1}=0 \;\textrm{ pour tout }
\X,\Y,\mathcal Z\in \Gamma(\T M).
$$

\subsection{Démonstration.}Soit $X+\xi$ une section de $\T
 M\lra M$. On rappelle qu'on identifie $\HH^\star$ aux formes linéaires sur  $T\ZZ(\QQ)$ qui sont  nulles sur $\V$. Notons  $\wh X+\wh \xi\in \Gamma(\HH\oplus\HH^\star)$ le relevé du champ
  $X+\xi$. Un tel champ est dit {\it basique}.

\quad\\
{\bf Proposition 1.} Soient  $A,B\in  \mathcal C^\infty(\V)$ deux champs de vecteurs
verticaux sur $\ZZ(\QQ)$. On se donne $X\in  \Gamma(TM)$ un champ de
vecteurs sur $M$ et $\xi \in  \Gamma(T^\star M)$ une 1-forme sur $M$.
\begin{enumerate}

\item $[A,B]\in \Gamma(\V)$

\item $[\wh X,A]\in \Gamma(\V)$

\item $[\wh X+\wh \xi ,\JJ A]=\JJ[\wh X+\wh \xi,A]$

\item $[\JJ(\wh X+\wh\xi),\JJ A]=\JJ[\JJ(\wh X+\wh \xi),A]$

\end{enumerate}

\quad\\
{\bf Preuve.}  Le premier point provient du fait que la
distribution verticale est l'espace tangent aux fibres de $\pi :
\ZZ(\QQ)\lra M$. Comme $\wh X$ est un champ relevé, le deuxième
point est immédiat. Et comme le transport parallèle suivant les
directions horizontales respecte l'orientation et la métrique sur
les fibres, elle respecte la structure complexe sur l'espace
tangent vertical, on a $[\wh X,\JJ A]=\JJ[\wh X,A]$. De plus, par
définition du crochet de Courant on vérifie que
$[\wh\xi,A]=0=[\wh\xi,\JJ A]$. Ce qui termine la preuve du
troisième  point. Le quatrième point s'obtient par "linéarité". En
effet soit $(\wh\X_1,\ldots,\wh\X_{8n})\in \HH\oplus\HH^\star$ une
base de champs  horizontaux basiques définie au-dessus de $\mathcal U$. Comme $\JJ$ stabilise
$\HH\oplus\HH^\star$, on note $[\JJ_{ij}]$ la matrice de la
restriction de $\JJ$ à $\HH\oplus\HH^\star$ dans cette base. En
utilisant le point 3, on a :
$$
\begin{array}{cccc}
&\JJ[\JJ\wh\X_j,A]&=&\JJ[\JJ_{ij}\wh\X_i,A]\\
&&=&\JJ\big(\JJ_{ij}[\wh\X_i,A]-A\wh\X_j\big)\\
&&=&\JJ_{ij}[\wh\X_i,\JJ A]-\JJ A\wh\X_j\\
\\\textrm{et}&[\JJ\wh\X_j,\JJ A]&=&[\JJ_{ij}\wh\X_i,\JJ A]\\
&&=&\JJ_{ij}[\wh\X_i,\JJ A]-\JJ A\wh\X_j.\; \square
\end{array}
$$

\quad\\
 {\bf Corollaire 1.} Le tenseur de Nijenhuis de $\JJ$ sur $\ZZ(\QQ)$ vérifie
 $\NN(\X,A)=0$ pour tout
$\X\in \Gamma(\HH\oplus\HH^\star)$ et pour tout $A\in \Gamma(\V)$.

\quad\\
{\bf Preuve.} Par linéarité, on peut supposer que $\X$ est un
champ basique. Le corollaire est alors une conséquence immédiate
des points 3 et 4 de la proposition 1. $\square$

\quad\\
{\bf Proposition 2.} Soient $\X, \Y\in \Gamma(\T M)$.
Au-dessus de notre ouvert $\mathcal U$,  la décomposition du champ
$[\wh \X,\wh \Y]$ en partie horizontale et verticale est donnée,
au point $u\in\ZZ(\QQ)$ par 
$$[\wh \X,\wh \Y] =\wh{[\X,\Y]}+[u,\RR(\X,\Y)].$$

\quad\\
{\bf Preuve.}  Soit $\mathcal P\lra M$ le $G$-fibré principal,
$\theta$ la $G$-connexion sur $\mathcal P$ et
 $ker\theta$ la distribution horizontale associée. Commençons par
 étudier le cas où $X$ et $Y$ sont deux champs de vecteurs sur $M$. Pour cela on note
  $\widetilde X,\widetilde Y$ leur relevé horizontal dans $\mathcal P$. La
décomposition du champ de vecteurs $[\widetilde X,\widetilde Y]$
en partie horizontale et verticale est donnée par (cf. \cite{O},
\cite{Bes87} chap 9) :
$$
[\widetilde X,\widetilde
Y]=\widetilde{[X,Y]}+(\theta\vert_{\V})^{-1}(\RR(X,Y)),$$
 où par définition, $(\theta\vert_{\V})^{-1}(\RR(X,Y))$ est le
 champ de vecteurs vertical sur $\mathcal P$ définie au point $p\in \mathcal P$ par :
 $$
\frac{d}{dt}\vert_{t=0}\Big(p.\,exp(t\RR(X,Y))\Big)=p.\,\RR(X,Y).
$$
D'autre part, la variété $\ZZ(\QQ)$ est le fibré associé de
fibre $\Sp^2$. Plus précisément, le groupe $G$ agit à droite sur
$\mathcal P\times\Sp^2$ :
$$
\begin{array}{ccc}
\mathcal P\times \Sp^2\times G&\lra&\mathcal P\times\Sp^2\\
(p,j,g)&\lms&(p.g,g^{-1}.j)=(p.g,gjg^{-1})
\end{array}
$$
et $\ZZ(\QQ)$ est le quotient de $\mathcal P\times\Sp^2$ par
cette action. On notera $\Pi$ la projection 
$$
\begin{array}{cccc}
\Pi &: \mathcal P\times\Sp^2&\lra&\ZZ(\QQ)\\
&(p,j)&\lms&u=p^{-1}jp.
\end{array}
$$
Comme au point $p\in\mathcal P$ on a
$d\Pi\Big(p.\RR(X,Y)\Big)=[u,\RR(X,Y)]$, au point $u\in\ZZ(\QQ)$
on a bien 
$$
[\wh X,\wh Y]=\wh{[X,Y]}+[u,\RR(X,Y)].
$$
Considérons maintenant $X$ un champ de vecteurs et $\eta$  une
1-forme sur $M$ alors par définition du crochet de Courant on voit
que $ [\wh X,\wh \eta]=\wh{[X,\eta]}$ et comme $\nabla_\eta=0$ on
a $\RR(X,\eta)=0$ soit :
$$
[\wh X,\wh \eta]=\wh{[X,\eta]}+[u,\RR(X,\eta)]. \square
$$

\quad\\
{\bf Corollaire 2.} Soient $U^\sharp\in \Gamma(\V^\star)$ une 1-forme
verticale et
 $\X\in  \Gamma(\HH\oplus \HH^\star)$ un champ de vecteurs et de formes
 horizontal. Au point $u\in\ZZ(\QQ)$, le tenseur de Nijenhuis $\NN(U^\sharp,\X)$
  est la 1-forme horizontale définie  pour tout $\overrightarrow{\Y}\in \Gamma(\HH)$ par
$$
 \NN(U^\sharp,\X)(\overrightarrow{\Y})=
 U^\sharp\left(\left[u,\RR\Big(\X\wedge\Y-u\X\wedge u\Y\Big)+
     u\RR\Big(u\X\wedge\Y+\X\wedge u\Y\Big)\right]\right).
 $$

 \quad\\
 {\bf Preuve.} Par définition du crochet de Courant on sait que
 $[U^\sharp,\X]=[U^\sharp,\overrightarrow{\X}]$ est une 1-forme.
 Soient $A\in \Gamma(\V)$ et $\overrightarrow{\Y}\in  \Gamma(\HH)$ deux
 champs de vecteurs,
au point $u\in\ZZ(\QQ)$ on a par définition 
$$
\begin{array}{ccc}
[U^\sharp,\X](A+\overrightarrow{\Y})&=&dU^\sharp(\overrightarrow{\X},\overrightarrow{\Y}+A)\\
&=&\overrightarrow{\X}.U^\sharp(A)-U^\sharp([\overrightarrow{\X},\overrightarrow{\Y}+A])\\
&=&\overrightarrow{\X}.U^\sharp(A)-U^\sharp([\overrightarrow{\X},A])
-U^\sharp([u,\RR(\X\w\Y)])
\end{array}
$$
Le point 3 de la proposition 1 nous assure alors que $[\JJ
U^\sharp,\X](A)=\JJ[U^\sharp,\X](A).$ La partie verticale de la
1-forme $\NN(U^\sharp,\X)$ est donc nulle.
 Pour la partie horizontale le calcul précédent
nous dit qu'au point $u$  on a
$$
\NN(U^\sharp,\X)(\overrightarrow{\Y})=U^\sharp\left(\left[u,\RR\Big(\X\wedge\Y-u\X\wedge
u\Y\Big)+
     u\RR\Big(u\X\wedge\Y+\X\wedge u\Y\Big)\right]\right).
  \;\square
 $$

\quad\\
{\bf Corollaire 3.} Soient $\X,\Y\in \Gamma(\HH\oplus\HH^\star)$ deux
champs horizontaux. Au point $u\in\ZZ(\QQ)$, le tenseur de
Nijenhuis $\NN(\X,\Y)$ est  le champ de vecteurs vertical défini
par 
$$\NN(\X,\Y)=-\left[u,\RR\Big(\X\wedge\Y-u\X\wedge u\Y\Big)+
     u\RR\Big(u\X\wedge\Y+\X\wedge u\Y\Big)\right].
$$

\quad\\
{\bf Preuve.}
  On note
$(\X_1,\ldots,\X_{8n})$ une  base orthonormée de $\T M$ définie
au-dessus de $\mathcal U$. La distribution $\mathcal
H\oplus\HH^\star$ est stable par $\JJ$, on notera $[\JJ_{ij}]$ sa
matrice dans la base relevée $(\wh{\X_1},\ldots,\wh{\X_{8n}})$.
Par définition du crochet de Courant 
$$\begin{array}{rll}
  \left[\mathbb J \widehat{\X_i},\mathbb
  J\widehat{\X_j}\right]&=&\overrightarrow{\JJ\wh{\X_i}}.(\JJ_{rj})\;\widehat{\X_r}
  -\overrightarrow{\JJ\wh{\X_j}}.(\JJ_{li})
\;\widehat{\X_l}+  \JJ_{li}\JJ_{rj}\left[\widehat{\X_l},\widehat{\X_r}\right]\\
&&-\JJ_{ri}d\JJ_{rj}+\JJ_{lj}d\JJ_{li}\\
  \,\left[\mathbb J\widehat{\X_i},\widehat{\X_j}\right]
  +\left[\widehat{\X_i},\mathbb
  J\widehat{\X_j}\right]  &=&-\overrightarrow{\wh{\X_j}}.(
  \JJ_{li})\;\widehat{\X_l}+\JJ_{li}\left[\widehat{\X_l},
  \widehat{\X_j}\right]+ \overrightarrow{\wh{\X_i}}.(
  \JJ_{rj})\;(\widehat{\X_r})
  +\JJ_{rj}\left[\widehat{\X_i},\widehat{\X_r}\right]\\
  &&+d\JJ_{ji}-d\JJ_{ij}.
\end{array}
$$
En utilisant la proposition 2, on en déduit qu'au point $u$, la
partie verticale de $\NN(\wh\X_i,\wh\X_j)$ vaut 
$$
-\left[u,\RR\Big(\X_i\wedge\X_j-u\X_i\wedge u\X_j\Big)+
     u\RR\Big(u\X_i\wedge\X_j+\X_i\wedge u\X_j\Big)\right].
$$
Pour la partie horizontale, on se donne $s$ une section locale de
$\ZZ(\QQ)\lra M$ telle que $s(m)=u$ et $(\nabla s)_m=0$. La
partie horizontale de $\NN(\wh\X_i,\wh\X_j)$ restreinte à la sous
variété $s(M)$ est égale au relevé horizontal du tenseur de
Nijenhuis $\NN(\X_i,\X_j)$ de $M$ munie de la structure presque
complexe généralisée induite par $s$. On se place dans des
coordonnées normales en $m$. En ce  point, on a alors :
$$\begin{array}{llc}<\NN(\X_i,\X_j),\X_k>&=&\TT(\X_i,\X_j,\X_k)-\TT(u\X_i,u\X_j, \X_k)\\
&& -\TT(u\X_i,\X_j,u\X_k)-\TT(\X_i,u\X_j,u\X_k)\Big).
 \end{array}$$
 Comme la connexion est sans torsion généralisée on en
déduit qu'au point $m$ on a $\NN(\X_i,\X_j)=0$. La partie
horizontale de $\NN(\wh\X_i,\wh\X_j)$ est nulle au point $(m,u)$
donc partout. $\square$

\quad\\
{\bf Preuve du théorème 1.} Comme les fibres de $\ZZ(\QQ)\lra M$
sont complexes, il est clair que quels que soient $A,B\in
 \Gamma(\V\oplus\V^\star)$ on a $\NN(A,B)=0$. Le théorème 1 est alors une
conséquence directe des corollaires 1, 2 et 3 et de la linéarité
de $\NN$. $\square$

\quad\\
{\bf Preuve du théorème 1 bis.} Par définition au point $(m,u)$ on
a :
$$
\begin{array}{ccr}
8\RR(\X^{1,0},\Y^{1,0})\ZZ^{1,0}&=& \Big(\RR( \X\w \Y-u\X\w u\Y)
-\RR( \X\w u\Y+u \X\w \Y)u\Big)\ZZ\\
&&-i\Big( \RR( \X\w u\Y+u \X\w \Y)+\RR( \X\w\Y-u\X\w u\Y)u\Big)\ZZ
\end{array}
$$
De sorte que la partie $(0,1)$ du vecteur
$8\RR(\X^{1,0},\Y^{1,0})\ZZ^{1,0}$ est nulle si et seulement si
$$\left[u,\RR\Big(\X\wedge
\Y-u\X\wedge u\Y\Big)+u\RR\Big(u\X\wedge\Y +\X\wedge
u\Y\Big)\right]=0.\; \square$$

\section{Variétés quaternioniques Kähler généralisées}
\subsection{Théorèmes d'intégrabilités}
 Soit $(M,g,\QQ)$ une 4n-variété riemannienne munie
d'une structure presque quaternionique hermitienne généralisée.
  En identifiant $\T M$ et $\T^\star M$ grâce à la pseudo-métrique
$<.,.>$, la métrique $g$ sur $\T M$ peut être vue comme un
endomorphisme
$G=\left(\begin{array}{cc}0&g^{-1}\\g&0\end{array}\right)$ de $\T
M=TM\oplus T^\star M$. Comme $G^2=Id$, on notera $C^\pm$ le
sous-espace propre de $G$ associé à la valeur propre $\pm1$. Si
au-dessus d'un ouvert $\mathcal U$, on se donne
$(\theta_1,\ldots,\theta_{4n})$ une base orthonormée de $TM$ et si
on note $(\theta_1^\star,\ldots,\theta_{4n}^\star)$ sa base duale,
alors 
$$
C^\pm=Vect\Big(\theta_1\pm \theta_1^\star,\ldots, \theta_{4n}\pm
\theta_{4n}^\star\Big).
$$
Comme les éléments de $\ZZ(\QQ)$ sont compatibles avec $g$, ils
stabilisent les espaces propres $C^\pm$. De plus comme
$\QQ-\{0\}\simeq \R^{+\star}\times\ZZ(\QQ)$, les éléments  de
$\QQ$ stabilisent aussi $C^\pm$. La matrice d'un élément
$u\in\QQ$ dans une base adaptée à $\T M=C^+\oplus C^-$ est donc de
la forme: $\left(\begin{array}{cc}
u^+&0\\0&u^-\end{array}\right)$. D'autre part, les projections
$p^\pm : C^\pm\lra TM$ sont des isomorphismes qui permettent de
descendre la distribution $\QQ$ en deux structures presque
quaternioniques sur $M$ que nous noterons
$Q^\pm:=p^\pm_\star(\QQ)$. Par abus de notation, nous écrirons
$u^\pm=p^\pm_\star(u)$.
 Enfin on définit l'application 
$$
\begin{array}{ccccccc}
f :&Q^-&\lra&\QQ&\lra&Q^+\\
&u^-&\lms&p^{-\star}(u^-)=\left(\begin{array}{cc}
u^+&0\\0&u^-\end{array}\right)&\lms&p^+_\star\circ
p^{-\star}(u^-)=u^+.
\end{array}
$$
{\bf Convention.} Nous dirons que $f : Q^-\lra Q^+$ est un
isomorphisme d'algèbre au sens où l'application $f$ se prolonge
naturellement un isomorphisme entre les algèbres $Q^-\oplus
Vect(Id)$ et $Q^+\oplus Vect(Id)$.

\quad\\
{\bf Corollaire 4.} La donnée d'une structure presque
quaternionique  hermitienne généralisée  sur $(M,g)$ revient à se
donner deux structures
 presque quaternioniques  hermitiennes $Q^+$ et $Q^-$ ainsi qu'un isomorphisme
 d'algèbre $f:Q^-\lra Q^+$. En particulier les espaces de twisteurs $\ZZ(\QQ)$, $Z(Q^+)$ et $Z(Q^-)$ sont difféomorphes :
$$
\xymatrix{& \ZZ(\QQ)\ar[ld]_{p^-_\star}\ar[rd]^{p^+_\star}&\\
    Z(Q^-) \ar[rr]^{f} & &Z(Q^+).
  }
  $$
\quad\\
{\bf Proposition 3.} Soit $\nabla$ la connexion de Levi-Civita de
$(M,g)$, $\QQ$ est une structure quaternionique Kähler généralisée
si et seulement si les deux conditions suivantes sont vérifiées:
\begin{enumerate}
\item  $Q^+$ et $Q^-$ sont des structures quaternioniques Kähler
sur $(M,g)$

\item $\nabla f(u^-)=f(\nabla u^-)$ pour tout $u^-\in Q^-$.
\end{enumerate}

\quad\\{\bf Remarque.} La première condition est automatiquement
vraie en dimension quatre.

\quad\\
{\bf Preuve.} La connexion de Levi-Civita stabilise $C^\pm$, donc
si elle stabilise $\QQ$, elle stabilise $Q^\pm$. Plus précisément
pour tout
$u=\left(\begin{array}{cc}u^+&0\\0&u^-\end{array}\right)\in\QQ$ on
a $\nabla u=\left(\begin{array}{cc}\nabla u^+&0\\0&\nabla
u^-\end{array}\right).$ De sorte que $\nabla
u\in\QQ\Longleftrightarrow
\left\{\begin{array}{l}\nabla u^+\in Q^+\\
\nabla u^-\in Q^-\\
\nabla u^+=f(\nabla u^-)\end{array}\right. .\quad\square$

\quad\\
Dans le cas particulier où $(M,g,\QQ)$ est une variété
riemannienne munit d'une structure
 quaternionique Kähler généralisée, on peut reformuler l'intégrabilité
 de la structure presque complexe généralisée
 $\JJ$ sur $\ZZ(\QQ)$ plus explicitement que dans le théorème 1.
 Comme dans le cas classique étudié par Salamon, l'intégrabilité de $\JJ$ est
 automatique dès lors que  $n>1$ :

\quad\\
{\bf Théorème 2.} Soit $(M,g,\QQ)$ une $4n$-variété quaternionique
Kähler généralisée avec $n>1$. La structure presque complexe
généralisée $\JJ$ sur son espace des twisteurs $\ZZ(\QQ)$ est
intégrable.

\quad\\
 Une 4-variété riemannienne orientée n'admet que
deux structures quaternioniques hermitiennes différentes : $\bw^+$ et $\bw^-$.
Quitte à renverser l'orientation, il y a donc deux classes de
structures quaternioniques Kähler généralisées suivant que
$Q^+=Q^-=\bw^+$ ou que $Q^+=\bw^+$ et  $Q^-=\bw^-$.

 \quad\\
 {\bf Théorème 3.} Soit $(M,g,\QQ)$ une 4-variété
 quaternionique Kähler généralisée telle que
  $Q^+=Q^-=\bw^+$. La
 structure presque complexe généralisée $\JJ$ sur $\ZZ(\QQ)$
 est intégrable si et seulement si l'une des conditions suivantes
 est vérifiée :
 \begin{enumerate}
 \item[a)] $R\vert_{\bw^+}=0$.

\item[b)] $f=Id$  et g anti-autoduale,

 \end{enumerate}

\quad\\{\bf Remarque.} Les seules 4-variétés riemanniennes
compactes qui vérifient la première contrainte sont les quotients
du tore plat ou d'une surface K3 \cite{Hit74}. D'autre part, si
$f=Id$ alors l'espace des twisteurs $\big(\ZZ(\QQ),\JJ\big)$
coïncident exactement avec l'espace des twisteurs
$\big(Z(\bw^+),\JJ\big)$ défini par Atiyah, Hitchin et Singer.

 \quad\\
 {\bf Théorème 4.} Soit $(M,g,\QQ)$ une 4-variété
 quaternionique  Kähler telle que
   $Q^+=\bw^+$ et $Q^-=\bw^-$.  La
 structure presque complexe généralisée $\JJ$ sur $\ZZ(\QQ)$
 est intégrable si et seulement si $g$ est localement conformément plate.

\quad\\
 Ces trois théorèmes permettent de retrouver
  l'intégrabilité de $\JJ$ pour les deux exemples donnés dans la section 2.

 \quad\\
 {\bf Exemple 1.} Si $(M,Q)$ est une variété quaternionique Kähler alors
 elle est naturellement quaternionique Kähler généralisée. Avec
 nos notations, cela revient à prendre
 $Q^+=Q^-=Q$ et  $f$ l'isomorphisme  identité. Pour $n>1$, l'intégrabilité de $\JJ$
  est automatique et on retrouve
le résultat de Salamon.  Pour $n=1$, l'intégrabilité de
 $\JJ$ équivaut à $g$ anti-autoduale
 et on retrouve le théorème d'Atiyah, Hitchin et Singer.

\quad\\
{\bf Exemple 2.} Soit $(M,g,I,J,K)$ une variété hyperkählérienne,
$ (\theta_1,\ldots,\theta_{4n})$ une base orthonormée locale de
$(TM,g)$ et $ (\theta_1^\star,\ldots,\theta_{4n}^\star)$ sa base
duale. Dans la base $(\theta_1,\theta_2,\ldots,\theta_{4n}^\star)$
de $\T M$ la structure quaternionique Kähler "tordue" décrite dans
la section 2.4, s'écrit
$$
\begin{array}{lll}
\QQ&=&\{a\J_I+b\J_w+c\J_{Iw}/(a,b,c)\in\R^3\}\\
&=&\left\{\left(\begin{array}{cc}aI&bJ+cK\\bJ+cK&aI\end{array}\right)/(a,b,c)\in\R^3\right\}.
\end{array}
$$
Dans la base
$(\theta_1+\theta_1^\star,\ldots,\theta_{4n}-\theta_{4n}^\star)$
de $C^+\oplus C^-$ cela donne donc:
$$
\QQ=\left\{\left(\begin{array}{cc}aI+bJ+cK&0\\0&aI-bJ-cK\end{array}\right)/(a,b,c)\in\R^3\right\}.
$$
Avec nos notations, cela correspond donc à prendre $Q^+=Q^-=Q$ et
$$
\begin{array}{lclc}
f:&Q^-&\lra&Q^+\\
&aI+bJ+cK&\lms&aI-bJ-cK
\end{array}
$$
Qu'on ait $n=1$ (comme la métrique est anti-autoduale et Ricci
plate) ou $n>1$, on rentre dans les hypothèses des théorèmes 2 et
3:  la structure presque complexe généralisée  $\JJ$ sur
$\ZZ(\QQ)$ est bien intégrable.
 De plus, on vérifie facilement que le type de $\JJ$ est
$2n+1$ au point $(m,\pm \mathcal J_I)$ et de type $1$ ailleurs.

\quad\\
La démonstration de ces trois théorèmes repose sur la proposition
suivante.

\quad\\
{\bf Proposition 4.} Soit $(M,g,\QQ)$ une variété quaternionique
Kähler généralisée et $R$ le tenseur de courbure de la connexion
de Levi-Civita. La structure presque complexe généralisée $\JJ$
sur $\ZZ(\QQ)$ est intégrable si et seulement si pour tous
champs de vecteurs $X,Y$ sur $M$ et pour tout élément $u\in
\ZZ(\QQ)$, les trois tenseurs suivants sont nuls :

\begin{enumerate}
\item[a)] $G_1(X,Y,u^+)=\left[u^+,R\Big(X\wedge Y
      -u^+X\wedge u^+Y\Big)+
        u^+R\Big(u^+X\wedge  Y+X\wedge u^+Y\Big)\right]. $

\item[b)] $G_2(X,Y,u^+)=\left[u^+,R\Big(X\wedge Y
      -u^+X\wedge u^-Y\Big)+
        u^+R\Big(u^+X\wedge  Y+X\wedge u^-Y\Big)\right]. $

\item[c)] $G_3(X,Y,u^-)=\left[u^-,R\Big(X\wedge Y
      -u^-X\wedge u^-Y\Big)+
        u^-R\Big(u^-X\wedge  Y+X\wedge u^-Y\Big)\right]. $
\end{enumerate}

\quad\\
{\bf Preuve.} La connexion de Levi-Civita $\nabla$ stabilise
$C^\pm$ donc pour tous $X,Y\in TM$, le tenseur de courbure
$\RR(X,Y)$ aussi. Plus exactement,  dans la base $\T M=C^+\oplus
C^-$, la matrice du tenseur de courbure $\RR(X,Y)$ s'écrit
$\left(\begin{array}{cc}R(X,Y)&0\\0&R(X,Y)\end{array}\right)$.
Au-dessus de l'ouvert $\mathcal U$ de $M$, on a la trivialisation
$$
\begin{array}{ccc}
\mathcal U\times\Sp^2&\lra&\pi^{-1}(\mathcal U)\subset \ZZ(\QQ)\\
(m,u^-)&\lms&
\left(m,\left(\begin{array}{ll}f(u^-)&0\\0&u^-\end{array}\right)\right).
\end{array}
$$
Les vecteurs verticaux sont donc de la forme
$\left(0,\left(\begin{array}{ll}f(A)&0\\0&A\end{array}\right)\right)$,
où $A$ est un vecteur tangent à $\Sp^2$. En utilisant cette
remarque et le théorème 1, on a que $\JJ$ est intégrable si et
seulement si le tenseur :
$$
G(\X,\Y,u^+)= \left[u^+,R\Big(\X\wedge\Y-u\X\wedge u\Y\Big)
+u^+R\Big(u\X\wedge\Y+\X\wedge u\Y\Big)\right]
$$
est  nul pour tous $\X,\Y\in\Gamma\Big(\T \ZZ(\QQ)\Big)$ et pour tout $u\in
\ZZ(\QQ)$. Suivant que  $\X$ et $\Y$ sont dans $C^+$ ou dans
$C^-$ on en déduit que l'annulation du tenseur $G$ est équivalente
à celle des tenseurs $G_1,G_2$ et $\tilde{G_3}$ avec:
$$
\tilde G_3(X,Y,u^+)=\left[u^+,R\Big(X\wedge Y
      -u^-X\wedge u^-Y\Big)+u^+R\Big(u^-X\wedge  Y+X\wedge u^-Y\Big)\right].
$$
Mais comme $\nabla u^+=f(\nabla u^-)$ et que
$[u,R(X,Y)]=-R(X,Y)u,\quad\forall u\in\bw^2$, on en déduit que
$[u^+,R(X,Y)]=f\Big([u^-,R(X,Y)]\Big)$ et donc que $\tilde
G_3(X,Y,u^+)=f\Big(G_3(X,Y,u^-)\Big) $. De sorte que $\tilde
G_3=0$ si et seulement si $G_3=0$. $\square$

\subsection{Démonstration du théorème 2.}  Pour démontrer le
théorème 2 nous utiliserons la proposition suivante:

\quad\\
{\bf Proposition 5.} Une 4n-variété quaternionique Kähler
généralisée avec $n>1$, vérifie nécessairement $f=Id$ ou $s=0$.

\quad\\
{\bf Preuve.} La démonstration de la proposition 5 repose sur le
lemme suivant:

\quad\\
{\bf Lemma 1~\cite{Bes87}.}  Soit $(M,g,Q)$ une 4n-variété
quaternionique Kähler avec $n>1$. On note $(I,J,K)$ une structure
presque hypercomplexe qui engendre $Q$ au-dessus d'un ouvert
$\mathcal U$ et $c=\displaystyle \frac{s}{2n(n+2)}$. Pour tous
champs de vecteurs $X,Y$ sur $M$ on a :
  \[
[I,R(X,Y)]=-cg(KX,Y)\,J+cg(JX,Y)\, K.
  \]
 Soit $(\II,\J,\KK)$ une structure presque
hypercomplexe généralisée  sur $(M,g)$ qui engendre localement
$\QQ$.  La proposition 3 nous assure que
$(\II^\pm,\J^\pm,\KK^\pm)$ sont deux structures quaternioniques
Kähler sur $(M,g)$. Le lemme 1 nous dit alors que:

$$
\begin{array}{ll}
 [\II^-,R(X,Y)]=-cg(\KK^-X,Y)\,\J^-+cg(\J^-X,Y)\, \KK^-\\

 [\II^+,R(X,Y)]=-cg(\KK^+X,Y)\,\J^++cg(\J^+X,Y)\, \KK^+.
\end{array}
 $$
En composant par $f$ la première égalité, on en déduit que :
$$
\left\{\begin{array}{lll}
 cg(\KK^-X,Y)&=& cg(\KK^+X,Y)\\
 cg(\J^-X,Y)&=& cg(\J^+X,Y)
\end{array}
\right.\Longleftrightarrow  f=Id \textrm{ ou } c=0.\quad \square
 $$
 Il est alors facile de démontrer le théorème 2. En effet
si $f=Id$, la variété $(\ZZ(\QQ),\JJ)$ coïncident avec la
variété $(Z(Q),\JJ)$ étudiée par Salamon, si bien que $\JJ$ est
une structure complexe intégrable sur $\ZZ(\QQ)$. D'autre part
dans le cas où $s=0$, il est clair en utilisant le lemme 1, que
$G_1=G_2=G_3=0$ et donc que $\JJ$ est intégrable.

\subsection{Démonstration du théorème 3.}
 On se place
maintenant en dimension quatre avec $Q^+=Q^-=\bw^+$.
 Le tenseur $G_1$ est la contrainte classique
d'intégrabilité qui apparaît dans Atiyah, Hitchin et Singer
\cite{AHS78}. Donc $G_1$ est nul  si et seulement si la
métrique $g$ est anti-autoduale. De même si $g$ anti-autoduale, on
a automatiquement $G_3=0$.  Il reste alors à étudier le tenseur
$G_2$,  sous l'hypothèse $g$ anti-autoduale. Pour cela on utilise
la proposition suivante qui est l'équivalent en dimension quatre
de la proposition 5.

\quad\\
{\bf Proposition 6.} Une 4-variété quaternionique Kähler
généralisée anti-autoduale  telle que $Q^+=Q^-=\bw^+$ vérifie
nécessairement  $f=Id$ ou $s=0$.

\quad\\
{\bf Preuve.} Sous l'hypothèse $g$ anti-autoduale, le tenseur de
courbure s'écrit:
$R=\left[\begin{array}{cc}\frac{s}{12}Id&\mathcal B\\
\mathcal B^\star&\mathcal W^-+\frac{s}{12}Id\end{array}\right]$.
De plus on sait que $[u^+,R(v)]=f\Big([u^-,R(v)]\Big)$ pour tout $
u\in \QQ$ et pour tout $v\in\bw^2$. Enfin, comme les éléments de $\bw^+$
commutent avec ceux de $\bw^-$, cela donne en particulier pour tout
$v\in\bw^+$:
$$
\begin{array}{llll}
&\displaystyle \frac{s}{12}[u^+,v]=\frac{s}{12}f\Big([u^-,v]\Big)&

 \Longleftrightarrow&\displaystyle \frac{s}{12}[u^+,v]=\frac{s}{12}[u^+,f(v)]\\

 &&\Longleftrightarrow& f=Id \textrm{ ou } s=0.\quad\square
 \end{array}
$$

\quad\\
 Si $f=Id$ alors $u^+=u^-$ et $G_1=G_2=G_3$. Donc $\JJ$
  est intégrable si et seulement si $g$ est
anti-autoduale. On considère donc le cas $f\neq Id$ et $\mathcal
W^+=s=0$.  Soit $(\theta_1,\ldots,\theta_4)$ une base orthonormée
au-dessus d'un ouvert $\mathcal U$ de $M$, cela nous fournit une
trivialisation locale de $\bw^+$ et de $\bw^-$ grâce aux deux
structures presque hypercomplexes suivantes
 $$
\left\{ \begin{array}{llr}
I^+&=&\theta_1\w\theta_2+\theta_3\w\theta_4\\
J^+&=&\theta_1\w\theta_3-\theta_2\w\theta_4\\
K^+&=&\theta_1\w\theta_4+\theta_2\w\theta_3
\end{array}
\right. \textrm{ et } \left\{ \begin{array}{llr}
I^-&=&\theta_1\w\theta_2-\theta_3\w\theta_4\\
J^-&=&\theta_1\w\theta_3+\theta_2\w\theta_4\\
K^-&=&-\theta_1\w\theta_4+\theta_2\w\theta_3
\end{array}
\right.
$$
On se place en un point $(m,u)$ de $\ZZ(\QQ)$ tel que $u^-=aI^++bJ^++cK^+$ et
$u^+=I^+$ alors:
$$\begin{array}{lll}
&G_2(\theta_1,\theta_1,u^+)=0\\
\Longleftrightarrow&\left[I^+, \mathcal
B\Big(-bK^--cJ^-\Big)+I^+\,
 \mathcal B\Big((a-1)I^-+bJ^--cK^-\Big)\right]=0&(1)\\
&\\
&G_2(\theta_2,\theta_2,u^+)=0\nonumber\\
\Longleftrightarrow&\left[I^+,\mathcal B\Big(bK^-+cJ^-\Big)+I^+\,
 \mathcal B\Big((a-1)I^--bJ^-+cK^-\Big)\right]=0&(2)\\
&\\
&G_2(\theta_3,\theta_1,u^+)=0\nonumber\\
\Longleftrightarrow&\left[I^+,\mathcal
B\Big((a-1)J^--bI^-\Big)+I^+\, \mathcal
B\Big((1-a)J^--cI^-\Big)\right]=0&(3)
\end{array}
$$
Comme $f\neq Id$ on a $a\neq 1$ et donc $(1)+(2)$ entraîne
$[I^+,\mathcal B(I^-)]=0$. Par symétrie, on en déduit que
$\mathcal
B=\left(\begin{array}{lll}x&0&0\\0&y&0\\0&0&z\end{array}\right)$
 est diagonale. Mais alors $(3)$ donne:
$$
\begin{array}{lc}
&[I^+,(a-1)yJ^++I^+(1-a)zK^+]=0 \Longleftrightarrow y=-z
\end{array}
$$
et par symétrie $z=-x$ et $x=-y$ soit $\mathcal B=0$.
Réciproquement, si $f$ quelconque, $g$ anti-autoduale et Ricci
plat alors l'image de $R$ est dans $\bw^-$. Comme les éléments de
$\bw^+$ commutent avec ceux de $\bw^-$ il est clair que
$G_1=G_2=G_3=0$ et donc que $\JJ$ est intégrable.

\subsection{Démonstration du théorème 4.} On se place toujours
en dimension quatre en supposant maintenant que $Q^+=\bw^+$ et
$Q^-=\bw^-$. Le tenseur $G_1$ est nul si et seulement si $g$ est
anti-autoduale et le tenseur $G_3$ est nul si et seulement si $g$
autoduale \cite{AHS78}. Il reste donc à étudier le tenseur $G_2$, sous
l'hypothèse $g$ localement conformément plate. Pour cela on
utilise la proposition et le lemme suivant:

\quad\\
{\bf Proposition 7.} Une 4-variété quaternionique Kähler
généralisée localement conformément plate telle que $Q^+=\bw^+$ et
$Q^-=\bw^-$ vérifie nécessairement $\mathcal
B=\displaystyle\frac{s}{12}f$.

\quad\\
{\bf Preuve.} Si la métrique $g$ est localement conformément
plate, la matrice du tenseur de courbure  est de la forme
$R=\left(\begin{array}{cc} \frac{s}{12}Id &\mathcal B\\\mathcal
B^\star&\frac{s}{12}Id\end{array}\right)$. Comme  pour tout $
u\in \QQ$ et pour tout $v\in\bw^2$ on a
$[u^+,R(v)]=f\Big([u^-,R(v)]\Big)$, en particulier pour tout $u^+\in\bw^+$ et pour tout $v\in\bw^-$ :
$$
\begin{array}{llll}
[u^+,\mathcal B(v)]=\displaystyle f\Big([u^-,\frac{s}{12}v]\Big)
 &\Longleftrightarrow&\displaystyle [u^+,\mathcal B(v)]=[u^+,\frac{s}{12}f(v)]\\
 &\Longleftrightarrow& \displaystyle \mathcal
B=\frac{s}{12}f.\quad \square
\end{array}
$$

\quad\\
{\bf Lemme 2.} Soit $(\II,\J,\KK)$ une structure presque
hypercomplexe généralisée qui engendre localement $\QQ$. En
choisissant convenablement la base
$(\theta_1,\theta_2,\theta_3,\theta_4)$, on peut supposer que:
$$
\left\{\begin{array}{l} I^\pm=\II^\pm\\
J^\pm=\J^\pm\\
K^\pm=\KK^\pm
\end{array}
\right.
$$
où $(I^\pm,J^\pm,K^\pm)$ est la base de $\bw^\pm$ définie dans la
section 4.3.

\quad\\
{\bf Preuve.}   Commençons par montrer qu'on peut choisir la base
$(\theta_1,\theta_2,\theta_3,\theta_4)$ de sorte que :
$$I^+=\II^+\textrm{ et } I^-=\II^-.
$$
 Les valeurs propres de l'endomorphisme $\II^+\II^-=\II^-\II^+$ sont  $-1$ et
$1$, on notera $E_{-1}$ et $E_1$  les espaces propres associés.  Soit $(\theta_1,\theta_2)$ une base orthonormée de $E_{-1}$
et $(\theta_3,\theta_4)$ une base orthonormée de $E_1$. Comme
$\II^+$ stabilise $E_{- 1}$ et est orthogonale, on a
$I^+\theta_1=\pm \theta_2$. Quitte à échanger $\theta_1$ et
$\theta_2$, on peut supposer que $\II^+\theta_1=\theta_2$. De même
on peut supposer que $\II^+\theta_3=\theta_4$. On a alors
$\II^+=\theta_1\w\theta_2+\theta_3\w\theta_4$. Mais comme
$\II^+\II^-\theta_1=-\theta_1$ et $\II^+\II^-\theta_3=\theta_3$ on
a aussi :
$$\II^-=\theta_1\w\theta_2-\theta_3\w\theta_4.$$
Maintenant en jouant sur le choix de $\theta_1\in E_{-1}$ nous
allons pouvoir conclure. Comme $f$ est un isomorphisme d'algèbre ,
sa matrice dans les bases $(I^-,J^-,K^-)$ de $\bw^-$ et
$(I^+,J^+,K^+)$ de $\bw^+$ associées à
$(\theta_1,\theta_2,\theta_3,\theta_4)$,  est de la forme :
$\left(\begin{array}{cc}1&\\&R_\alpha\end{array}\right)$ où
$R_{\alpha}$ est  la rotation d'angle $\alpha$. Notons
$R_{-\alpha}$ la rotation d'angle $-\alpha$ dans le plan $E_{-1}$.
On vérifie directement que la base
$\Big(R_{-\alpha}(\theta_1),R_{-\alpha}(\theta_2),\theta_3,\theta_4\Big)$
convient. $\square$

\quad\\
On se place dans la base $(\theta_1,\theta_2,\theta_3,\theta_4)$
donnée par le lemme 2. La proposition 6 nous dit que la matrice du
tenseur de courbure dans la base $(I^+,J^+,K^+,I^-,J^-,K^-)$ est
$R=\displaystyle\frac{s}{12}\left(\begin{array}{cc}Id&Id\\Id&Id\end{array}\right)$.
Pour finir la preuve du théorème 4, il suffit de vérifier à la
main que  $G_2(\theta_i,\theta_j,I^+)=0$, pour tous $\theta_i,\theta_j$.

\section{Nouveaux exemples}
 Pour illustrer ces théorèmes, on peut donner de
nouveaux exemples de structures quaternioniques Kähler
généralisées. L'idée est de partir de deux structures
quaternioniques Kähler $Q^+, Q^-$ sur une 4n-variété riemannienne
$(M,g)$ et d'un isomorphisme d'algèbre $f:Q^-\lra Q^+$. On peut
alors définir la structure presque quaternionique hermitienne généralisée
$\QQ\subset \OO(\T M)$ dont l'écriture dans $\T M=C^+\oplus
C^-$ est
$$
\QQ_f=\left\{\left(\begin{array}{cc}f(u)&0\\0&u\end{array}\right)/u\in
Q^-\right\}.
$$
 Si on a bien choisi nos structures de
départ et notre isomorphisme, on obtient alors une structure
quaternionique Kähler
 généralisée sur $(M,g)$ et  une structure presque complexe généralisée $\JJ$ sur
 $\ZZ(\QQ_f)$ intégrable.

\quad\\
{\bf Application 1.} Soit  $(M,g,I,J,K)$  une variété hyperkählérienne
et $Q=\{aI+bJ+cK/(a,b,c)\in\R^3\}$.  On oriente $Q$ en décrétant
que $I,J,K$ est directe. Soit $f$ n'importe quelle isométrie
directe de $Q$, la proposition 3 nous assure que $\QQ_f$ est une
structure quaternionique Kähler généralisée sur $(M,g)$ et on a
automatiquement l'intégrabilité de la structure complexe $\JJ$ par
les théorèmes 2 et 3. C'est une généralisation de l'exemple 2. En
particulier on peut voir l'exemple 2 comme une déformation de la
structure canonique et la terminologie structure quaternionique
Kähler généralisée "tordue" prend tout son sens.

\quad\\
{\bf Application 2.} Plus généralement étant donnée $(M,g,Q)$ une
variété quaternionique Kähler munit d'une structure complexe
kählérienne $I\in Q$, considérons $f_\theta : Q\lra Q$ la rotation
 d'angle $\theta$ autour de l'axe $(-I,I)$.

\quad\\
{\bf Lemme 3.}  La distribution $\QQ_{f_\theta}$ est une structure
 quaternionique Kähler généralisée sur $(M,g)$.

\quad\\
{\bf Preuve.} Il suffit là encore d'utiliser la proposition 3.
Soit $(I,J,K)$ une structure presque hypercomplexe qui engendre
localement $Q$. Comme $I$ est kählérienne, il existe une 1-forme
$a$ sur $M$ telle que $ \left\{\begin{array}{l}
 \nabla I=0\\
 \nabla J=aK\\
 \nabla K=-aJ
 \end{array}
 \right.
 $
et on a bien $\nabla f_\theta(u)=f_\theta(\nabla u)$ pour tout
$u\in Q$. $\square$

\quad\\
 Ainsi, si $M$ est une surface d'Enriques ou une surface hyperelliptique \cite{BHPVdV} alors le théorème 3
 nous assure que la
structure $\JJ$ sur $\ZZ(\QQ_{f_\theta})$ est intégrable. On a construit une famille de déformation de
structures complexes généralisées sur la variété $\ZZ(\QQ_f)\simeq Z(\bw ^+)$ dont les types valent:
$$
\begin{array}{lcc}
a)& 2n+1 &\textrm{ si }\theta\in 2\pi\mathbb Z\\

b)& \left\{\begin{array}{cl} 2n+1 &\textrm{ aux points } (m, \J_{\pm I})\\
1 & \textrm{ partout ailleurs}\end{array}\right.&\textrm{ si
}\theta\notin 2\pi\mathbb Z
\end{array}
$$
 Pour donner d'autres exemples, nous utiliserons le lemme
suivant.

\quad\\
{\bf Lemme 4.} Soient $M$ est une 4-variété parallélisable,
$(\theta_1,\theta_2,\theta_3,\theta_4)$ une base orthonormée de
champs de vecteurs et $\Gamma_{ij}^k$ les symboles de Christoffel
de la connexion de Levi-Civita:
$$
\nabla_{\theta_i}\theta_j=\sum_k\Gamma_{ij}^k\theta_k.
$$
On conserve la notation de la partie 4.3, où la base
$(\theta_1,\theta_2,\theta_3,\theta_4)$ de $TM$ définit
 une base $(I^\pm,J^\pm,K^\pm)$ de $\bw^\pm$.  Pour tout $i$, on
a
$$
\begin{array}{lclll}
\nabla_{\theta_i}I^+&=&
(\Gamma_{i1}^4+\Gamma_{i2}^3)J^+&+&(-\Gamma_{i1}^3+\Gamma_{i2}^4)K^+\\

\nabla_{\theta_i}J^+&=&
(\Gamma_{i3}^2-\Gamma_{i1}^4)I^+&+&(\Gamma_{i3}^4+\Gamma_{i1}^2)K^+\\

\nabla_{\theta_i}K^+&=&
(\Gamma_{i4}^2+\Gamma_{i1}^3)I^+&+&(\Gamma_{i4}^3-\Gamma_{i1}^2)J^+\\

\\

\nabla_{\theta_i}I^-&=&
(-\Gamma_{i1}^4+\Gamma_{i2}^3)J^-&+&(-\Gamma_{i1}^3-\Gamma_{i2}^4)K^-\\

\nabla_{\theta_i}J^-&=&
(\Gamma_{i3}^2+\Gamma_{i1}^4)I^-&+&(-\Gamma_{i3}^4+\Gamma_{i1}^2)K^-\\

\nabla_{\theta_i}K^-&=&
(-\Gamma_{i4}^2+\Gamma_{i1}^3)I^-&+&(-\Gamma_{i4}^3-\Gamma_{i1}^2)J^-\\
\end{array}
$$

\quad\\
 {\bf Preuve.} Par définition on a:
 \[
\begin{array}{cll}
(\nabla_{\theta_i}I^+)\theta_1&=&\nabla_{\theta_i}(I^+\theta_1)-I^+\nabla_{\theta_i}\theta_1\\
&=&\nabla_{\theta_i}\theta_2-I^+(\Gamma_{i1}^k\theta_k)\\
&=&\left[\begin{array}{l} \Gamma_{i2}^1\\
\Gamma_{i2}^2\\
\Gamma_{i2}^3\\
\Gamma_{i2}^4\end{array}\right]
-\left[\begin{array}{r} -\Gamma_{i1}^2\\
\Gamma_{i1}^1\\
-\Gamma_{i1}^4\\
\Gamma_{i1}^3\end{array}\right]=
\left[\begin{array}{c} 0\\
0\\
\Gamma_{i2}^3+\Gamma_{i1}^4\\
\Gamma_{i2}^4-\Gamma_{i1}^3
\end{array}\right]
\\
\Longleftrightarrow
\nabla_{\theta_i}I^+&=&(\Gamma_{i2}^3+\Gamma_{i1}^4)J^++(\Gamma_{i2}^4-\Gamma_{i1}^3)K^+.
\end{array}
\]
Les autres expressions se démontrent de la même manière. $\square$

\quad\\
{\bf Application 3.} Les surfaces hyperelliptiques sont
parallélisables \cite{Des1}, on se fixe une métrique riemannienne
$g$ et une base orthonormée $(\theta_1,\theta_2,\theta_3,\theta_4)$. On peut
alors définir la structure presque quaternionique généralisée sur
$(M,g)$ telle que $Q^\pm=\bw^\pm$ et associée à l'isomorphisme
$f_\theta$ dont la matrice dans les bases $(I^-,J^-,K^-)$ et
$(I^+,J^+,K^+)$
 est
$f_\theta=\left(\begin{array}{ccc}1&0&0\\0&\cos\theta&-\sin\theta\\
0&\sin\theta&\cos\theta\end{array}\right)$.

\quad\\
{\bf Lemme 5.} Pour toutes surfaces hyperelliptiques, on peut
choisir une métrique $g$ et une base orthonormée
$(\theta_1,\theta_2,\theta_3,\theta_4)$ de sorte que, quel que soit
$\theta$, la distribution $\QQ_{f_\theta}$ soit une structure
quaternionique Kähler généralisée. De plus la structure
presque complexe généralisée $\JJ$ sur son espace des twisteurs
généralisées $\ZZ(\QQ_{f_\theta})$ est intégrable, son type est constant égale à deux.

\quad\\
{\bf Preuve.}  Commençons par rappeler la classification des
surfaces hyperelliptiques (\cite{BHPVdV} Chap.V.5). Pour cela,
notons $I_0$ la structure complexe canonique sur $\R^2$ qui permet
d'identifier $\R^2$ et $\C$. Soit $\T_1$ le tore $\C/\Z+i\Z$ et
$\T_2$ le tore $\C/\Gamma$, o\`u $\Gamma$ est un r\'eseau de $\C$.
On note $z_1$ et $z_2$ les coordonn\'ees complexes respectivement
sur $\T_1$ et sur $\T_2$.  \`A diff\'eomorphisme pr\`es il y a
sept types de
 surface hyperelliptique. Dans tous les cas, c'est le quotient de
 $\T_1\times \T_2$ par
un des groupes  d'automorphismes $G$ suivant :
\newpage
\quad\\
\begin{tabular}{c|c|c|c}
Type&$\Gamma$&$G$& Action de $G$ sur $\T_1\times\T_2$\\
\hline
 &&&\\
 1& $\Z+i\Z$&$\Z/2\Z$&$g_1(z_1,z_2)=(z_1+\frac{1}{2},-z_2)$\\
 &&&\\
2& $\Z+i\Z$&$\Z/2\Z\oplus\Z/2\Z$&$\left\{\begin{array}{l}
g_1(z_1,z_2)=(z_1+\frac{1}{2},-z_2)\\
g_2(z_1,z_2)=(z_1+\frac{i}{2},z_2+e_1)\;{\rm avec}\;2e_1=0
\end{array}\right.$\\
&&&\\
\hline &&&\\
3 &$\Z+j\Z$&$\Z/3\Z$&$g_1(z_1,z_2)=(z_1+\frac{1}{3},jz_2)$\\
&&&\\

 4&$\Z+j\Z$&$\Z/3\Z\oplus\Z/3\Z$&$\left\{\begin{array}{l}
g_1(z_1,z_2)=(z_1+\frac{1}{3},jz_2)\\
g_2(z_1,z_2)=(z_1+\frac{i}{3},z_2+e_1)\;{\rm avec}\;je_1=e_1
\end{array}\right.$\\
&&&\\

\hline &&&\\
5&$\Z+i\Z$&$\Z/4\Z$&$g_1(z_1,z_2)=(z_1+\frac{1}{4},iz_2)$\\
 &&&\\

6&$\Z+i\Z$&$\Z/4\Z\oplus\Z/2\Z$&$\left\{\begin{array}{l}
g_1(z_1,z_2)=(z_1+\frac{1}{4},iz_2)\\
g_2(z_1,z_2)=(z_1+\frac{i}{2},z_2+e_1)\;{\rm avec}\;ie_1=e_1
\end{array}\right.$\\
&&&\\
\hline &&&\\
7&$\Z+j\Z$&$\Z/6\Z$&$g_1(z_1,z_2)=(z_1-\frac{1}{6},-jz_2)$\\

\end{tabular}

\quad\\
\quad\\
Quel que soit le type de la surface hyperelliptique $M$, la
métrique canonique sur $\C^2$ descend en une métrique riemannienne
$g$ sur $M$.  Nous allons maintenant construire une base
orthonormée de champs de vecteurs. Le premier tore $\T_1$ est le
quotient $\C/\Z+i\Z$. On note $(e_1,e_2)$ la base canonique de
$\C\simeq\R+i\R$. Le deuxi\`eme tore est le quotient
$\T_2=\C/\Gamma$. On note ici $(e_3,e_4)$ la base canonique de
$\C\simeq\R+i\R$ et $R_\theta$ la rotation d'angle $\theta$ dans
le plan $\C\simeq Vect(e_3,e_4)$ orienté par $e_3\w e_4$. On
d\'efinit alors sur $\C\oplus\C$ les quatre champs de vecteurs :
\[
\left\{
\begin{array}{l}
\theta_1(z_1,z_2)=R_{2\pi\,Re\,z_1}(e_3)\\
\theta_2(z_1,z_2)=R_{2\pi\,Re\,z_1}(e_4)\\
\theta_3(z_1,z_2)=e_1\\
\theta_4(z_1,z_2)=e_2
\end{array}
\right.
\]
Ces champs passent naturellement au quotient en quatre champs sur
$\T_1\times\T_2$, quelque soit le r\'eseau $\Gamma$ d\'efinissant
 $\T_2$. De plus quelque soit le type de la surface
hyperelliptique, ces champs passent au quotient sur
$\T_1\times\T_2/G$. On a ainsi construit une base de champs de
vecteurs pour chaque surface hyperelliptique.

 Quel que soit le type
de la surface hyperelliptique, les structures complexes
kählériennes $I^\pm=I_0\oplus\pm I_0$ sur $\R^2\oplus \R^2$ passe
au quotient. Les structures $I^\pm$ sont donc kähleriennes et
comme $\nabla_{\theta_i}\theta_4=0=\nabla_{\theta_i}\theta_3$, le
lemme 4 donne :
$$\begin{array}{cc}
\left\{\begin{array}{ccc}
\nabla_{\theta_i}I^+&=& 0\\

\nabla_{\theta_i}J^+&=&\Gamma_{i1}^2K^+\\

\nabla_{\theta_i}K^+&=&-\Gamma_{i1}^2J^+

\end{array}
\right. \quad \left\{\begin{array}{ccc}

\nabla_{\theta_i}I^-&=&0\\

\nabla_{\theta_i}J^-&=&\Gamma_{i1}^2K^-\\

\nabla_{\theta_i}K^-&=&-\Gamma_{i1}^2J^-
\end{array}
\right.
\end{array}
$$
La proposition 3 et le théorème 4 nous assure alors que
$\QQ_{f_\theta}$ et $\JJ$ sont intégrables, pour tout $\theta$.
$\square$

\quad\\
{\bf Application 4.} Considérons ici une variété localement
conformément plate $M$ obtenue par produit de $\Sp^1$ avec une
surface de dimension 3 à courbure sectionnelle constante non
nulle, c'est-à-dire de la forme $\Sp^3\times \Sp^1$ ou
$\Sigma\times \Sp^1$, avec $\Sigma$ une 3-variété hyperbolique.
Soit $(\theta_1,\theta_2,\theta_3)$ une base orthonormée globale
de $\Sp^3$ ou de $\Sigma$ et $\theta_4$ un champ de vecteur
unitaire sur $\Sp^1$. Comme dans la partie 4.3, on note encore
$(I^\pm,J^\pm,K^\pm)$ la base de $\bw^\pm$ associé à la base
$(\theta_1,\theta_2,\theta_3,\theta_4)$ et on considère $\QQ_f$ la
distribution définie par $Q^\pm=\bw^\pm$ et
$f(aI^-+bJ^-+cK^-)=aI^++bJ^++cK^+$. Comme $\Gamma_{i4}^j=0$ pour
tout $i,j$, le lemme 4 et la proposition 3 nous assure que $\QQ_f$
est une structure quaternionique Kähler généralisée et le théorème
4 nous assure l'intégrabilité de la structure complexe généralisée
$\JJ$ sur $\ZZ(\QQ)$. Son type est constant égale à 2.

\quad\\
{\bf Application 5.}  Le tore plat $\T^4$ admet deux structures
hyperkählériennes : $\bw^+$ et $\bw^-$.  Le tore plat $\T^8$,
produit de deux tores $\T^4$, admet donc quatre structures
hyperkählériennes produits différentes. Soient $Q^+$ et $Q^-$ deux
telles structures et $f$ n'importe quel isomorphisme  entre $Q^-$
et $Q^+$. La distribution $\QQ_f$ est alors une structure
quaternionique Kähler généralisée et la structure presque complexe
généralisée $\JJ$ sur $\ZZ(\QQ_f)$ est intégrable.

\quad\\
{\bf Remarque.} Dans ces trois derniers exemples, on a construit une structure complexe généralisée sur $\ZZ(\QQ)$. Mais  cette variété admet également deux structures complexes différentes : celle donnée par $Z(Q^+)$ et celle donnée par $Z(Q^-)$. Cela nous fournit donc des exemples de variétés bi-hermitiennes qui sont de plus complexes généralisées.

\quad\\

\quad\\
Guillaume DESCHAMPS,\\
\quad\\
 Université de Brest, UMR 6205,\\ Laboratoire de
Mathématiques de
Bretagne Atlantique\\
6 avenue Victor le Gorgeu\\
CS 93837,\\
29238 Brest cedex 3\\
(France)


\begin{thebibliography}{10}


\bibitem{AHS78}
M.F. Atiyah, N.J. Hitchin, I.M. Singer.
\newblock {\em Self-duality in four-dimensional {R}iemannian
geometry.}
\newblock {Proc. Roy. Soc. London Ser. A} {\bf 362}, 425--461 (1978)


\bibitem{Ale}
D.V. Alekseevskii.
\newblock {\em Riemanniann manifolds with exceptional holonomy groups.}
\newblock {Functional Anal. Appl.} {\bf 2}, 97-105 (1968)



\bibitem{BHPVdV}
W.  Barth, K. Hulek, C.  Peters,
              A. Van de Ven.
\newblock {\em Compact complex surfaces.}
\newblock {Ergebnisse der Mathematik und ihrer Grenzgebiete. 3. Folge. A
              Series of Modern Surveys in Mathematics},
\newblock Springer-Verlag, Berlin (2004)

\bibitem{Ber}
M. Berger.
\newblock {\em Remarques sur le groupe d'holonomie des variétés riemanniennes.}
\newblock {C. R. Acad. Sci. Paris} {\bf 262}, 1316-1318 (1966)


\bibitem{Bes87}
A.L. Besse.
\newblock {\em Einstein manifolds.}  Ergebnisse der Mathematik
  und ihrer Grenzgebiete (3),
\newblock Springer-Verlag, Berlin (1987)



\bibitem{Bre}
A. Bredthauer.
\newblock {\em Generalized Hyperkähler geometry and supersymmetry.}
\newblock   Nucl. Phys. B {\bf 773}, 172-183 (2007)




\bibitem{Cou}
T.J. Courant.
\newblock {\em Dirac manifolds.}
\newblock   Trans. Amer. Math. Soc. {\bf 319}, 631-661 (1990)

\bibitem{Des1}
G. Deschamps.
\newblock {\em  Espaces twistoriels et structures complexes non standards.}
\newblock Publicacions Math. {\bf 52}, 435-457 (2008)


\bibitem{Des2}
G. Deschamps.
\newblock {\em  Espace des twisteurs des structures complexes généralisées.}
\newblock Mathematische Zeitschrift {\bf 279}, 703-721 (2015)


\bibitem{DM1}
J. Davidov, O. Mushkarov.
\newblock {\em  Twistor spaces of generalized complex structures.}
\newblock J. Geom. Phys. {\bf 56}, 1623-1636 (2006)



\bibitem{DM2}
J. Davidov, O. Mushkarov.
\newblock {\em  Twistorial construction of generalized Kähler manifolds.}
\newblock  J. Geom. Phys. {\bf 57}, 889-901 (2007)

\bibitem{GS}
R. Glover, J. Sawon.
\newblock{\em Generalized twistor spaces for hyperkähler manifolds.}
\newblock J. London Math. Soc. {\bf 91}, 321-342 (2015)

\bibitem{Gua2}
M. Gualtieri.
\newblock {\em Generalized complex geometry.}
\newblock Ph.D. thesis, St John's college, University of Oxford,
 arXiv: math/0401221,   107 pages (2003)


\bibitem{Gua3}
M. Gualtieri.
\newblock {\em Branes on Poisson varieties.}
\newblock The many Facets of Geometry, Oxford Scholarship Online Monographs, 368-395 (2010)

\bibitem{Hit74} N.~J. Hitchin.
\newblock {\em Compact four-dimensional Einstein manifolds.}
\newblock Journal of Differential Geometry {\bf 9}, 435-441
(1974)



\bibitem{Hit}
N.J. Hitchin.
\newblock {\em Generalized Calabi-Yau manifolds.}
\newblock  Q. J. Math. {\bf 54}, 281-308 (2003)










\bibitem{NN}
A. Newlander, L. Nirenberg.
\newblock {\em Complex analytic coordinates in almost complex manifolds.}
\newblock Ann. of Math. {\bf 65}, 391--404 (1957)


\bibitem{O}
B.  O'Neil.
\newblock {\em The fundamental equations of a submersion.}
\newblock  Michigan Math. J. {\bf 13}, 459-469 (1966)








\bibitem{Pan}
R. Pantilie.
\newblock {\em  Generalized quaternionic manifolds.}
\newblock   Ann. Mat. Pura ed Appl.
{\bf 193}, 633-641, (2014)







\bibitem{Pen76}
R. Penrose.
\newblock {\em  Nonlinear gravitons and curved twistor
theory.}
\newblock General Relativity and Gravitation {\bf 7}, 31-52 (1976)







\bibitem{Sal}
S. Salamon.
\newblock {\em Quaternionic Kähler manifolds.}
\newblock  Invent. Math. {\bf 67}, 143-171 (1982)


\bibitem{Sal1}
S. Salamon.
\newblock {\em Quaternionic manifolds.}
\newblock  Symposia Matematica {\bf 26}, 139-151 (1982)







\bibitem{ST}
I.M. Singer, J.A. Thorpe.
\newblock {\em The curvature of {$4$}-dimensional Einstein spaces.}
\newblock Global Analysis, paper in honor of K. Kodaira, Princeton University Press, 355-365 (1969)









\end{thebibliography}
\end{document}